\documentclass[11pt, reqno]{amsart}

\DeclareMathAccent{\mathring}{\mathalpha}{operators}{"17}

 \usepackage{color}

\newcommand{\mysection}[1]{\section{#1}
      \setcounter{equation}{0}}

\newtheorem{theorem}{Theorem}[section]
\newtheorem{lemma}[theorem]{Lemma}

\newtheorem{corollary}[theorem]{Corollary}

\theoremstyle{definition}
\newtheorem{assumption}{Assumption}[section]

\theoremstyle{remark}
\newtheorem{remark}{Remark}[section]

\newtheorem{example}{Example}[section]

\newcommand\cbrk{\text{$]$\kern-.15em$]$}}
\newcommand\opar{\text{\raise.2ex\hbox{${\scriptstyle | }$}\kern-.34em$($} }
\newcommand{\tr}{\text{\rm tr}\,}

 \makeatletter
 \def\dashint{%
 \operatorname%
 {\,\,\text{\bf--}\kern-.98em\DOTSI\intop\ilimits@\!\!}}
 \makeatother

 \newcommand{\WO}{\overset{\rm o}{ W}\,\!}

\newcommand\bR{\mathbb{R}}

\newcommand\bL{\mathbb{L}}
\def\bM{\mathbb{M}}

\newcommand\cS{\mathcal{S}}

\newcommand\dist{{\rm dist}\,}

\begin{document}

\title[Existence theorem for
equations with relaxed convexity]
{On the existence of $W^{2}_{p}$ 
solutions for fully nonlinear elliptic
equations under  relaxed convexity assumptions}

\author{N.V. Krylov}
\thanks{The  author was partially supported by
an NSF Grant}
\email{krylov@math.umn.edu}
\address{127 Vincent Hall, University of Minnesota,
 Minneapolis, MN, 55455}



\begin{abstract}
We establish the existence and uniqueness of solutions
of fully nonlinear elliptic second-order equations
like $H(v,Dv,D^{2}v,x)=0$
in smooth domains without requiring $H$
to be convex or concave with respect to the second-order derivatives.
Apart from ellipticity nothing is required of $H$ at points at which
 $|D^{2}v|\leq K$,
where $K$ is any given constant. For large $|D^{2}v|$ some kind
of relaxed convexity assumption with respect to $D^{2}v$ mixed with a VMO condition
with respect to $x$ are still imposed. The solutions are sought
in Sobolev classes.
\end{abstract}

\maketitle

\mysection{Introduction}
                                   \label{section 12.14.1}

In the literature,   interior $W^2_p, p>d$, a priori
estimates for a class of
fully nonlinear uniformly elliptic equations in $\bR^{d}$ of the form
\begin{equation}
                                                  \label{11.30.4}
 H(v,Dv,D^2v,x)=0
\end{equation}
were first obtained by Caffarelli in \cite{Caf89} (see also
\cite{CC95}).
 Adapting his technique, similar interior 
a priori estimates were proved by Wang
\cite{Wa92} for parabolic equations. In the same paper, a boundary
estimate is stated but without a proof; see Theorem 5.8 there.  By
exploiting a weak reverse H\"older's inequality, the result of
\cite{Caf89} was sharpened by Escauriaza in \cite{Es93}, who obtained the
interior $W^2_p$-estimate for the same equations allowing
$p>d-\varepsilon$, with a small constant
$\varepsilon$ depending only on the ellipticity constant and $d$. Quite
recently, Winter \cite{Wi09} further extended this technique to establish
the corresponding boundary a priori estimates as well as the {\em
$W^2_p$-solvability\/} of the associated boundary-value problem. It is
also worth noting that a solvability theorem in the space
$W^{1,2}_{p,\text{loc}}(Q)\cap C(\bar Q)$ can be found in \cite{CKS00} for
the boundary-value problem for fully nonlinear parabolic equations. 
The above mentioned results of \cite{CKS00} and \cite{Wi09}
are proved under the assumption that $H$ is convex in $D^{2}v$
and in all papers mentioned above
 a small oscillation assumption in the integral sense is
imposed on the operators; see, for instance,
\cite[Theorem 1]{Caf89}. However, as pointed out in
\cite[Remark 2.3]{Wi09} and in \cite{Kr10}
(see also \cite[Example 8.3]{CKS00} for a relevant discussion),
this assumption turns out to be equivalent to a small oscillation condition in the $L_\infty$ sense,
which, particularly in the
{\em linear\/} case, is the same as what is required in the
classical $L_{p}$-theory based on the Calder\'on--Zygmund
theorem when one first investigates the case of
constant coefficients and then by using perturbation
method and partitions of unity passes to the case that
the coefficients are uniformly sufficiently close to
continuous ones. Thus, it seems to the author that the results in
\cite{Caf89,Wa92,Es93,CKS00,Wi09} mentioned above are in general not
formally applicable
to the operators under our Assumption
\ref{assump1} on the ``main'' part $F$ of $H$,
 in which local oscillations are measured in a certain average sense 
allowing rather rough discontinuities.
It is still possible that the {\em methods\/} developed in the above
cited articles can be used to obtain our results. In our opinion,
our method, which is quite different from
theirs,
is somewhat simpler and leads to the results faster.

So far \cite{DKL}  is
 the only  article  where fully nonlinear elliptic
and parabolic equations in smooth {\em domains\/} with VMO ``coefficients''
were shown to be solvable in (global) Sobolev classes. There the a priori
estimates are obtained under assumptions which are stronger than
ours but yet very similar. However, the solvability is
proved only under the assumption that $H$ {\em is\/} convex in $
D^{2}v$.
Here we prove the conjecture stated in \cite{DKL} that the convexity
assumption and a kind of bounded inhomogeneity assumption
can be dropped in the case of elliptic equations.
The author intends to do the same for parabolic equations in a subsequent
paper.

The results obtained in this article generalize and  contain the Sobolev
space theory of {\em linear\/} equations
 with VMO coefficients, which was developed about twenty years ago   by
Chiarenza, Frasca, and Longo in  \cite{CFL1,CFL2} for non-divergence form
elliptic equations, and later in \cite{BC93} by Bramanti and Cerutti for
parabolic equations. The proofs in these references are based on
explicit representations of second-order derivatives
through certain singular integrals, on  the
Calder\'on--Zygmund theorem and the Coifman--Rochberg--Weiss commutator
theorem. For further related results, we refer the reader to the book
\cite{MaPaSo00} and reference therein.
There are no explicit solutions for nonlinear equations and
this method cannot be applied. We use a different approach
described in \cite{DKL}  and \cite{Kr10}.

  As in \cite{DKL} we assume that $H$ is represented as the sum
of two functions: main part $F(D^{2}v,x)$ and a subordinated part
$G(v,Dv,D^{2}v,x)$. In \cite{DKL} the function $G$ is supposed to grow
sublinearly with respect to $|D^{2}v|$,
so that, as far as a priori estimates are concerned, one need only 
  estimate the $W^{2}_{p}$-norm of $v$ through the $L_{p}$-norm
of $F(D^{2}v,x)$. 

There are two differences between our
assumptions about $F(D^{2}v,x)$ and the ones made in \cite{DKL}.
First we do not suppose that $F$ is convex in $D^{2}v$. Actually,
the convexity of $F$ was never used in \cite{DKL} either. More important 
difference is that we do not assume that $F$ is positive homogeneous
with respect to $D^{2}v$. Instead, we assume that $F(0,x)=0$
(which holds automatically in \cite{DKL}) and that, in a sense, $F$ can be
approximated by convex functions for large $|D^{2}v|$.
One the one hand, we get an a priori estimate (in elliptic case)
under weaker conditions than in \cite{DKL} and, on the other hand,
without this generalization we would not be able to prove
the existence of solutions.

Relaxing the assumptions on $F$ leads to impossibility
of employing the usual localization techniques
 even in the case that $F$ is independent of $x$
 when, for instance, we want to estimate $I:=F(D^{2}(\zeta v))
-\zeta F(D^{2}v)$, where $\zeta\geq0$ is a smooth cut-off
function. It was used in \cite{DKL} that $\zeta F(D^{2}v)=
F(\zeta D^{2}v)$ and then the Lipschitz continuity of
$F$ guaranteed an estimate of $I$ through lower order terms.
Therefore, a new argument was needed, actually, avoiding using
partitions of unity and localizations altogether. 

After we obtain the necessary a priori estimates
we derive our existence theorem from a very general existence
theorem in which there is no assumptions on the structure of
$H$ or on its convexity (see Theorem \ref{theorem 9.23.1}).
Instead, equation  \eqref{11.30.4} is modified by using a parameter
$K\geq0$, so that when $K\to\infty$ the modified equation
becomes \eqref{11.30.4}. The reader may be surprised by the fact
that the assertions of Theorem \ref{theorem 9.23.1} 
apart from   estimates \eqref{1.13.2} and \eqref{2.28.1}
are proved in a basic case in \cite{Kr13}
with the help of finite-difference approximations without 
using any results from the theory of linear or fully nonlinear
elliptic equations. It is exactly the structure
of the approximating equations what prevented us
from using positive homogeneous $F$.

\mysection{Main result}

In this article, we consider elliptic  equations
\begin{equation}
                                                \label{7.29.1}
H[v](x):= H(v( x),D v( x),D^{2}v( x), x)=0
\end{equation}
in subdomains of $\bR^{d } $, where
$$
\bR^{d}=\{x=(x_{1},...,x_{d}):x_{1},...,x_{d}\in
\bR=(-\infty,\infty)\}.
$$
In \eqref{7.29.1}
$$
D^{2}u=(D_{ij}u),\quad Du=(D_{i}u),\quad
D_{i}=\frac{\partial}{\partial x_{i}},\quad D_{ij}=D_{i}D_{j}.
$$
We introduce $\cS$ as the set of symmetric $d\times d$ matrices,
fix a constant 
$\delta\in(0,1]$, and  set
$$
\cS_{\delta}=\{a\in\cS:\delta|\xi|^{2}\leq a_{ij}\xi_{i}\xi_{j}
\leq\delta^{-1}|\xi|^{2},\quad \forall \xi\in\bR^{d}\},
$$
where and everywhere in the article the summation convention is enforced
unless specifically stated otherwise.

Recall that   Lipschitz
continuous functions are almost everywhere differentiable.

\begin{assumption}
                                    \label{assumption 9.23.1}
 The function $H(u,x)$,
$u=(u',u'')$, 
$$
u'=(u'_{0},u'_{1},...,u'_{d})
\in\bR^{d+1},\quad u''\in\cS,
$$
 is measurable with respect to $x$ for any $u$,
and Lipschitz continuous
in $u$ for every $x\in\bR^{d}$. For any $x$, at all 
points of differentiability of $H(u,x)$ with respect to $u$
$$
(H_{u''_{ij}}(u,x))\in \cS_{\delta},\quad
|H_{u'_{k}}(u,x)|\leq K_{0},\quad k=1,...,d,
\quad 0\leq-H_{u'_{0}}(u,x)\leq K_{0}.
$$
where $K_{0}$ is a fixed constant.
\end{assumption}

Next we assume that there are two
functions $F(u,x)=F(u'',x)$ and $G(u,x)$ such that
$$
H =F +G .
$$
 
\begin{example}
                                          \label{example 11.29.1}
One can take
$F(u'',x)=H(0,u'',x)$ and $G=H-F$. 
Since we will require later that $F(0,x)=0$, one can then
take $F(u'',x)=H(0,u'',x)-H(0,x)$ and $G=H-F$. However, we are
not bound by these choices.
\end{example}

\begin{assumption}
                                    \label{assumption 10.5.1}

The function $G(u',u'' ,x)$, $u''\in\cS,u'\in\bR^{d+1}$,
 is nonincreasing
in $u'_{0}$ and  
$$
|G(u',u'',x)|\leq  K_{0}|u'|+\bar{G}(x).
$$

\end{assumption}
\begin{remark}
                                        \label{remark 11.7.1}
In \cite{DKL} a less restrictive assumption
is imposed on  $G$ allowing it to grow sublinearly
with respect to $u''$. However, this part of $G$
from \cite{DKL} can be absorbed into $F$
on account of increasing $t_{0}$ in our
Assumption \ref{assump1}. 

\end{remark}

Observe that, owing to Assumption \ref{assumption 9.23.1},
Assumption \ref{assumption 10.5.1} is satisfied
in Example \ref{example 11.29.1}
for the first choice of $F$ and $G$
and  $\bar{G}\equiv0 $.

In contrast with \cite{DKL} we do not 
suppose that $F$ is positive homogeneous of degree
one with respect to $u''$. This generalization is,
actually, necessary in order for 
 the method we prove our main result to go through.
However, we have to pay for that by having a more
complicated
 VMO (vanishing mean oscillation)  assumption containing
 a constant $\theta\in (0,1]$ to be specified later.
We also do not assume that $F$ is convex in $u''$.
The combination of
Assumption \ref{assumption 10.5.1} and the following
one we call  ``a relaxed convexity
assumption on $H$ and a VMO condition on $F$''.

 Set
$$
B_{r}(x)=\{y\in\bR^{d}:|x-y|<r\},\quad B_{r}=B_{r}(0)
$$
and for Borel $\Gamma\subset \bR^{d}$ denote by $|\Gamma|$
the volume of $\Gamma$.
Let $\Omega$ be an open bounded subset of $\bR^{d}$
with $C^{2}$ boundary.

\begin{assumption}
                                \label{assump1}
There exist  $R_0\in(0,1]$ and $t_{0}\in[0,\infty)$
such that for any
$r\in (0, R_0]$  and  $z\in \Omega$
one can find a {\em convex\/} function $\bar{F} (u'' )=
\bar{F}_{z,r} (u'' )$ (independent
of $x $)  such that

(i) We have $\bar{F}(0)=0$ and at all points of differentiability
of $\bar{F}$ we have $(\bar{F}_{u''_{ij}})\in\cS_{\delta}$;

(ii)
For any $u''\in\cS$ with $|u''|=1$ we have
\begin{equation}
                                                \label{7.30.2}
\int_{\Omega\cap B_{r}(z)}\sup_{t>t_{0}}t^{-1}
|F(tu'' ,x)-\bar{F}(tu'')| \,dx\leq \theta
|\Omega\cap B_{r}(z)|,
\end{equation}
where for $u''\in\cS$ by $|u''|$ we mean $ \tr ^{1/2}(u''u'')  $;

(iii) The function $F$ is Lipschitz continuous with respect to $u''$
with Lipschitz constant $K_{0}$, measurable
with respect to $x$, and $F(0,x)\equiv0$.
\end{assumption}

Here is our main result.

\begin{theorem}
                                    \label{theorem 10.5.1}
Let $p>d$ and assume
 that $\bar G\in L_{p}(\Omega)$.
Then there exists a constant $\theta\in (0,1]$, depending
only on $d$, $p$, $\delta$, and  $  \Omega$,
such that, if Assumption \ref{assump1} is satisfied with this
 $\theta$, then 

(i) For any $g\in W^{ 2}_p(\Omega)$ there exists a unique
$u\in W^{ 2}_p(\Omega)$ satisfying \eqref{7.29.1} and such
that $u-g\in \WO^{ 2}_p(\Omega)$.

(ii) We have
\begin{equation}
                                \label{eq16.01}
\|u\|_{W^{2}_p(\Omega)}\le N
\| \bar G\|_{L_p(\Omega)}+N\|g\|_{W^{1,2}_p(\Omega)}
 +N t_{0} ,
\end{equation}
where $N$ depends only on $K_{0}$, $d$, $p$, $\delta$,  
$R_{0}$, and
  $ \Omega$.
\end{theorem}

Here $W^{2}_p(\Omega)$ denotes the set of
functions $v$ defined on
$\Omega$ such that $v$, $Dv
 $, and $D^2v $ 
are in
$L_p(\Omega)$, and $\WO^{ 2}_p(\Omega)$ is the set of all functions
$v\in
W^{ 2}_p(\Omega)$ such that $v$ vanishes on $\partial \Omega$.

To prove   the uniqueness part of the theorem
introduce $\bL_{\delta,K_{0}}$ as the collection of operators
$$
Lu=a^{ij}D_{ij}u+b^{i}D_{i}u-cu
$$
 with measurable coefficients such that at all points
$a=(a^{ij})\in\cS_{\delta}$, $|b^{i}|\leq K_{0}$,
$i=1,...,d$,
$0\leq c\leq K_{0}$. 

It is a well-known fact that
owing to Assumption \ref{assumption 9.23.1}
for any $u,v\in W^{2}_{p}(\Omega)$ there exists
an operator $L\in\bL_{\delta,K_{0}}$ such that  
$H[u]-H[w]=L(u-w)$. Then uniqueness in Theorem
\ref{theorem 10.5.1} follows from the Alexandrov maximum  principle. 

The remaining assertions
of the  theorem are proved in Section \ref{section 11.29.1}
after we prove necessary a priori estimates.

\begin{remark}
                                         \label{remark 12.8.2}
The parameter $\theta$ in Theorem \ref{theorem 10.5.1}
depends on $p$ and we cannot guarantee that it stays
bounded away from zero for all $p>d$. Our arguments
are only valid if we take $\theta$ sufficiently small
and as $p\to\infty$, $\theta$ should go to zero.
\end{remark}

\begin{remark}
                                         \label{remark 12.8.1}
For a Borel set $\Gamma\subset\bR^{d}$ with nonzero 
Lebesgue measure and locally summable $f$ denote
$$
\dashint_{\Gamma}f(x)\,dx=\frac{1}{|\Gamma|}
\int_{\Gamma}f(x)\,dx,
$$
where $|\Gamma|$ is the volume of $\Gamma$. Then for $z\in\Omega$
and $r>0$ introduce
$$
\hat{F}(u'')
=\hat{F}_{z,r}(u'')=\dashint_{\Omega\cap B_{r}(z)}F(u'',x)\,dx.
$$

Observe that
$$
|F(tu'',x)-\hat{F}(tu'')|\leq |F(tu'',x)-\bar{F}(tu'')|
+\dashint_{\Omega\cap B_{r}(z)}|\bar{F}(tu'')-F(tu'',y)|\,dy,
$$
which implies that
\begin{equation}
                                             \label{12.8.2}
\dashint_{\Omega\cap B_{r}(z)}
\sup_{t>t_{0}}t^{-1}
|F(tu'' ,x)-\hat{F}(tu'')| \,dx\leq 2\theta.
\end{equation}
Thus, one can be tempted to always take $\hat{F}$ as $\bar{F}$.
However, there is no guarantee that $\hat{F}(u'')$ is convex
in $u''$.

\end{remark}

\begin{remark}
                                         \label{remark 11.29.1}
Under the above assumptions the function
$H(0,u'',x)-H(0,x)$ does not necessarily satisfy
Assumption \ref{assump1} (with $H(0,u'',x)-H(0,x)$ in place of
$F(u'',x)$), so that this choice of $F$ and $G$ in Example
\ref{example 11.29.1} may not be optimal. The simplest example
in case $d=2$
is given by 
$$
H(0,u'',x)=G(x)\wedge|u''_{11}|+2u''_{11}+u''_{22},
$$
where $G(x)= x_{1}  ^{-\alpha}$ for $x_{1}>0$ and
$G(x)= 0$ for $x_{1}\leq0$
with a small $\alpha>0$, so that $G$ is summable to a high power.

Indeed,  assume that $0\in\Omega$. Then
for $z=0$, small $r>0$,
 and $u''_{11}=1$ the left-hand side of \eqref{12.8.2}
(with $H(0,u'',x)-H(0,x)$ in place of
$F(u'',x)$)
becomes
$$
\dashint_{ B_{r} }
\sup_{t>t_{0}} 
\big| 1\wedge(G(x)/t)-\dashint_{ B_{r} }1\wedge(G(y)/t)\,dy\big| \,dx
$$
$$
\geq\dashint_{ B_{r} }
\big| 1\wedge(G(x)/t_{0})-\dashint_{ B_{r} }1
\wedge(G(y)/t_{0})\,dy\big| \,dx,
$$
which for $r\leq t_{0}^{1/\alpha}$ equals
$$
\dashint_{ B_{r} }
 | I_{x_{1}>0}-\dashint_{ B_{r} }I_{y_{1}>0}\,dy | \,dx
=\dashint_{ B_{r} }
 | I_{x_{1}>0}-1/2 | \,dx=1/2.
$$
Hence, \eqref{7.30.2} cannot be satisfied with small $\theta$
and
the natural choice for $F$ in this example is $2u''_{11}+u''_{22}$.

\end{remark}
\begin{remark}
                                         \label{remark 11.29.2}
In condition \eqref{7.30.2} no restriction is imposed
on $F(u'',x)$ for $|u''|\leq t_{0}$. But even for large
$|u''|$
the function $F$  satisfying  Assumption \ref{assump1}
need not be even locally convex. An example can be constructed
looking at the case that $d=2$ and
  $F(u'',x)=2u''_{11}+u''_{22}+f(|u''_{11}| )$, where
 $f$ is any sublinearly growing function with $|f'|\leq1$
and such that $f(0)=0$.

\end{remark}

In addition to the examples presented in Remarks \ref{remark 11.29.1}
and \ref{remark 11.29.2}
we give one more.

\begin{example}
                                             \label{example 12.7.1}
Let $A$ and $B$ be some countable sets and
assume that for $\alpha\in A$, $\beta\in B$, and $x\in\bR^{d}$
we are given functions $a^{\alpha }(x)$, $b^{\alpha\beta}(x)$,
$c^{\alpha\beta}(x)$, and $f^{\alpha\beta}(x)$ with values in
$\cS_{\delta}$, $\bR^{d}$, $[0,\infty)$, and $\bR$,
respectively. Assume that these functions are measurable
in $x$,  $b^{\alpha\beta}$ and $c^{\alpha\beta}$ 
are bounded, and
$$
\bar{G}:=\sup_{\alpha,\beta}|f^{\alpha\beta}|\in L_{p}(\Omega).
$$

Consider the following Isaacs equation
\begin{equation}
                                                \label{12.9.1}
H(v,Dv,D^{2}v,x)=0,
\end{equation}
where
$$
H(u,x):=
\inf_{\beta\in B}\sup_{\alpha\in A}
\big[\sum_{i,j=1}^{d}
a^{\alpha }_{ij}(x)u''_{ij}+\sum_{i =1}^{d}b^{\alpha\beta}_{i}(x)u'_{i}
-c^{\alpha\beta}(x)u'_{0}+f^{\alpha\beta}(x)\big].
$$
Our measurability, boundedness, and countability assumptions guarantee that
$H$ is measurable in  $x$ and Lipschitz continuous in $u$. One can
also easily check that 
at all points of differentiability $(H_{u''_{ij}})\in\cS_{\delta}$.
Next assume that there is an $R_{0}\in(0,\infty)$ such that
for any point $z\in\Omega$ and $r\in(0,R_{0}]$ one can find
$\bar{a}^{\alpha}\in\cS_{\delta}$ (independent of $x$) such that
$$
\sup_{\alpha\in A}\dashint_{\Omega\cap B_{r}(z)}
|a^{\alpha}(x)-\bar{a}^{\alpha}|\,dx\leq\theta,
$$
where $\theta$ is taken from Theorem \ref{theorem 10.5.1}.

Then we claim that the assertions (i) and (ii) of 
Theorem \ref{theorem 10.5.1}  
hold true and estimate \eqref{eq16.01} holds with $t_{0}=0$.

To prove the claim introduce
$$
F(u'',x)=\sup_{\alpha\in A}
 \sum_{i,j=1}^{d}
a^{\alpha }_{ij}(x)u''_{ij},\quad G=H-F.
$$
Notice that Assumption \ref{assump1} is satisfied 
with $t_{0}=0$ and  
$$
\bar{F}(u''):=\sup_{\alpha\in A}
 \sum_{i,j=1}^{d}\bar
a^{\alpha\ }_{ij} u''_{ij}
$$
because these functions are convex, positive homogeneous
of degree one
with respect to $u''$ and, for $|u''|=1$,
$$
\dashint_{\Omega\cap B_{r}(z)}|
F(u'',x)-\bar{F}(u'')|\,dx
\leq \sup_{\alpha\in A}\dashint_{\Omega\cap B_{r}(z)}\big|
\sum_{i,j=1}^{d}
[a^{\alpha }_{ij}(x)-\bar{a}^{\alpha}]u''_{ij}\big|
$$
$$
\leq\sup_{\alpha\in A}\dashint_{\Omega\cap B_{r}(z)}
|a^{\alpha}(x)-\bar{a}^{\alpha} |\,dx\leq\theta.
$$

On can easily check that Assumption \ref{assumption 10.5.1}
is satisfied as well and this proves our claim. 
\end{example}

In the proofs of various results in this article we use
the symbol $N$ sometimes with indices to denote constants
which may change from one occurrence to another and
we do not always specify on which data these  constants
depend. In these cases the reader should remember
that, if in the statement of a result there are constants
called $N$ which are claimed to depend only on certain
parameters, then in the proof of the result
the constants $N$ also depend only on the same
parameters unless specifically stated otherwise.

\mysection{Interior a priori estimates for the simplest
equation}
                                    \label{section 10.7.1}

In this and the following section we assume that
$F(u'')$ is a convex function of $u''$ (independent of $x$)
such that $F(0)=0$ and at all points of differentiability
of $F$ we have $(F_{u''_{ij}})\in\cS_{\delta}$.

\begin{lemma}
                                        \label{lemma 10.5.1}
  There exists an
$\alpha=\alpha(d,\delta)\in(0,1)$ such that for any 
 $\phi\in C(\partial B_{2})$ there exists a unique
$v\in C(\bar{B}_{2})\cap C^{2+\alpha}_{loc}(B_{2})$ satisfying
$$
F(D^{2}v)=0\quad\text{in}\quad B_{2},\quad v=\phi\quad
\text{on}\quad\partial B_{2}.
$$
Furthermore,
$$
|D^{2}v(x)-D^{2}v(y)|\leq N|x-y|^{\alpha}
\sup_{\partial B_{2}}|\phi|
$$
as long as $x,y\in B_{1}$,
where $N$ depends only on $\delta$ and $d$.
\end{lemma}

This lemma is a somewhat weaker version of
 Theorem 4.1 in \cite{Sa94}. Even though
the author of \cite{Sa94} attributes this lemma to
Evans-Krylov (see \cite{Ev}, \cite{Kr_82}) and it can be indeed extracted from the results
of chapter 5 of \cite{Kr85}, in the above clear and convenient
 form it is
 stated and proved in \cite{Sa94}.   
In what follows by $\alpha$ we mean the constant in
Lemma \ref{lemma 10.5.1} until further notification.

\begin{lemma}
                                          \label{lemma 12.16.01}

Let $r\in(0,\infty)$, $\nu\geq2$ and let $\phi\in
C(\partial B_{\nu r})$. Then there exists
a unique $v\in C(\bar{B}_{\nu r})\cap C^{2+\alpha}_{loc}
(B_{\nu r})$ such that
$$
F(D^{2}v)=0\quad\text{in}\quad B_{\nu r},\quad v=\phi\quad
\text{on}\quad\partial B_{\nu r}.
$$
Furthermore,
$$
\dashint_{B_{r}}\dashint_{B_{r}}|D^{2}v(x)-D^{2}v(y)|
\,dxdy\leq N(d,\delta)\nu^{-2-\alpha} r^{-2}\sup_{\partial
B_{\nu r}}|\phi|.
$$
\end{lemma}

Proof. Dilations show that it suffices to concentrate on
$r=2/\nu$. In that case the existence of solution
follows from Lemma \ref{lemma 10.5.1}, which also implies that
for $x,y\in B_{2/\nu}\subset B_{1}$
$$
|D^{2}v(x)-D^{2}v(y)|\leq N\nu^{-\alpha}\sup_{B_{2}}|\phi|.
$$
It only remains to observe that
$$
\dashint_{B_{2/\nu}}\dashint_{B_{2/\nu}}|D^{2}v(x)-D^{2}v(y)|
\,dxdy\leq\sup_{x,y\in B_{2/\nu}}|D^{2}v(x)-D^{2}v(y)|.
$$
  The lemma is proved.

The following is a slight
 generalization of the main  result of
\cite{FHL} proved in case $u=0$ on $\partial B_{1}$.
Lemma \ref{lemma 12.16.2} follows from
  Theorems 1.8 and 2.2
of \cite{Kr12.2} when $\gamma=\gamma_{0}$.
For arbitrary $\gamma\in(0,\gamma_{0}]$
one obtains the result by using H\"older's inequality.

\begin{lemma}
                                          \label{lemma 12.16.2}

There are   constants $\gamma_{0}\in(0,1]$ and $N$,
depending only on $\delta$, $K_{0}$,  and $d$, such that 
for any $L
\in\bL_{\delta,K_{0}}$, $\gamma\in(0,\gamma_{0}]$,  and $u\in 
W^{2}_{d,loc}(B_{1})\cap C(\bar B_{1})  $
 we have
$$
 \dashint_{B_{1}}(|D^{2}u|^{\gamma}+|Du|^{\gamma})
\,dx \leq N\big(\dashint_{B_{1}}|Lu|^{d}\,dxdt\big)^{\gamma/d}
 +N\sup_{\partial B_{1}}|u|^{\gamma} .
$$
\end{lemma}

Below in this section
 by $\gamma_{0}$ we always mean the constant in Lemma
\ref{lemma 12.16.2}. By using dilations we get the following.

\begin{corollary}
                                     \label{corollary 10.5.2}
For any $r\in(0,\infty)$, $L\in\bL_{\delta,0}$,
$\gamma\in(0,\gamma_{0}]$, and $u$ belonging to
$W^{2}_{d,loc}(B_{r})\cap C(\bar B_{r})  $
 we have
$$
 \dashint_{B_{r}}(|D^{2}u|^{\gamma}+r^{-\gamma}|Du|^{\gamma})
\,dx \leq N\big(\dashint_{B_{r}}|Lu |^{d}\,dxdt\big)^{\gamma/d}
 +Nr^{-2\gamma}\sup_{\partial B_{r}}|u|^{\gamma} .
$$

\end{corollary}

We keep following the example of notation given in
\eqref{7.29.1} and set
$$
F[u](x)=F( D^{2}u(x) ).
$$

\begin{corollary}
                                     \label{corollary 10.7.1}
There exist  $\gamma\in(0,\gamma_{0}]$ and $N$
depending only on $\delta$, $K_{0}$, $d$, and $\Omega$ such that
for any  $L\in\bL_{\delta,K_{0}}$, and $u\in 
W^{2}_{d,loc}(\Omega)\cap C(\bar \Omega)  $
 we have
\begin{equation}
                                       \label{9.5.06}
 \int_{\Omega}(|D^{2}u|^{\gamma}+|Du|^{\gamma})
\,dx \leq N\|Lu\|_{L_{d}(\Omega)}^{\gamma}
 +N\sup_{\partial \Omega}|u|^{\gamma} .
\end{equation}

\end{corollary}
 Indeed, one can represent $\bar{\Omega}$ as a finite union
of the closures of $C^{2}$-domains each of which
admits a one-to-one $C^{2}$ mapping on $B_{1}$
with $C^{2}$ inverse. Then after changing coordinates
one can use Lemma \ref{lemma 12.16.2} applied to appropriately 
changed operator $L$. For the transformed operator the
constants $\delta$ and $K_{0}$
may change but still will only depend
on $\delta, K_{0},d$, and $\Omega$. Then after combining the results
of application of Corollary \ref{corollary 10.5.2} one
obtains \eqref{9.5.06} with $\bar{\Omega}$ in place of
$\partial\Omega$. However, Alexandrov's estimate
shows that this replacement can be avoided on account
of, perhaps, increasing the first $N$ on the right
in \eqref{9.5.06}.

In what follows by $\gamma$ we mean the constant from Corollary
\ref{corollary 10.7.1}.

 For $\rho>0$ introduce
$$
\Omega_{\rho}=\{x:\rho(x) >\rho\},
$$
where
$$
\rho(x)=\dist(x,\Omega^{c}).
$$
\begin{lemma}
                                      \label{lemma 10.5.2}
 
Let $r\in(0,\infty)$ and $\nu\in(2,\infty)$. Then for any
$u\in W^{2}_{d,loc}(\Omega)$ and  $z\in\Omega_{\nu r}$ we have
$$
\big(\dashint_{B_{r}(z)}\dashint_{B_{r}(z)}
|D^{2}u(x)-D^{2}u(z)|^{\gamma}\,dxdy\big)^{1/\gamma}
$$
\begin{equation}
                                           \label{10.6.1}
\leq N\nu^{d/\gamma}
\big(\dashint_{B_{\nu r}(z)}
|F[u]|^{d}\,dx\big)^{1/d}+
N\nu^{-\alpha}\big(\dashint_{B_{\nu r}(z)}|D^{2}u|^{d}\,dx
\big)^{1/d},
\end{equation}
where $N$ depends only on $d$ and $\delta$.
\end{lemma}

Proof. Take a point $z\in\Omega_{\nu\rho}$ 
 and define
$v$ to be a unique $C(\bar{B}_{\nu r}(z))
\cap C^{2+\alpha}_{loc}(B_{\nu r}(z))$-solution of
equation $F[v]=0$ in $B_{\nu r}(z)$ with boundary condition
$v=u$ on $\partial B_{\nu r}(z)$. Such a function exists
by Lemma \ref{lemma 12.16.01} applied after shifting the origin.
Furthermore, $v(x)-b^{i}x_{i}-c$ satisfies the same equation
for any constant $b^{i},c$. Hence by Lemma \ref{lemma 12.16.01}
and H\"older's inequality
$$
I_{r}(z):=\big(
\dashint_{B_{r}(z)}\dashint_{B_{r}(z)}|D^{2}v(x)-D^{2}v(y)|^{\gamma}
\,dxdy\big)^{1/\gamma}
$$
$$
\leq N \nu^{-2-\alpha} r^{-2}\sup_{x\in\partial
B_{\nu r}(z)}|u(x)-(D_{i}u)_{B_{\nu r}}x_{i}-u_{B_{\nu r}}|.
$$
By Poincar\'e's inequality (see, for instance,
Lemma 2.1 in \cite{Kr10}) the last supremum is dominated
by a constant times
$$
\nu^{2} r^{2}\big(\dashint_{B_{\nu r}(z)}|D^{2}u|^{d}\,dx
\big)^{1/d}.
$$
It follows that
\begin{equation}
                                            \label{10.6.3}
I_{r}(z)\leq N\nu^{-\alpha}\big(\dashint_{B_{\nu r}(z)}|D^{2}u|^{d}\,dx
\big)^{1/d}.
\end{equation}

Next, the function $w=u-v$ is of class $C(\bar{B}_{\nu r}(z))
\cap W^{2}_{p,loc}(B_{\nu r}(z))$ and for an operator  
 $L\in\bL_{\delta,0}$ we have $F[u]-F[v]=L(u-v)$, 
$ L(u-v)=F[u]$. Moreover, $w=0$ on $\partial B_{\nu r}$.
Therefore,
by Corollary \ref{corollary 10.5.2}
$$
\dashint_{B_{r}(z)}|D^{2}w|^{\gamma}\,dx
\leq \nu^{d}\dashint_{B_{\nu r}(z)}|D^{2}w|^{\gamma}\,dx
$$
$$
\leq N\nu^{d}\big(\dashint_{B_{\nu r}(z)}
|F[u]|^{d}\,dx\big)^{\gamma/d}.
$$
Upon combining this with \eqref{10.6.3} we come 
to \eqref{10.6.1} and the lemma is proved.

\mysection{Boundary a priori estimates in the simplest case}
                                    \label{section 10.7.2}

We suppose that the assumptions stated in the beginning of
Section \ref{section 10.7.1} are satisfied
and set
$$
\bR^{d}_{+}=\{x\in\bR^{d}:x=(x_{1},x'), x_{1}>0\},\quad
B_{r}^{+}=\{|x|<r:x_{1}>0\}.
$$
 
For real numbers $z\geq 0$ denote
$$
B_{r}^{+}(z)=\{x\in\bR^{d}_{+} : |x-z e_{1}|<r\},
$$
where $e_{1}$ is the first basis vector in $\bR^{d}$.

\begin{lemma}
                                                     \label{osca}
There exists an
$\alpha=\alpha(d,\delta)\in(0,1)$ such that if
 $z,r>0$, 
$\nu\geq  16 $,
$$
u\in C(\bar{B}^{+}_{\nu r}(z))\cap\bigcap_{\rho<\nu r}
W^{2}_{d}(B^{+}_{\rho}(z)),
$$
and $u$ vanishes for $x_{1}=0$, then 
we have
$$
 \big(\dashint_{B^{+}_{r}(z)}\dashint_{B^{+}_{r}(z)}
|D^{2}u(x)-D^{2}u(y)|^{\gamma}\,dxdy\big)^{1/\gamma}
$$
\begin{equation}
                                              \label{10.8.4}
\leq N\nu^{d/\gamma}
\big(\dashint_{B^{+}_{\nu r}(z)}
|F[u]|^{d}\,dx\big)^{1/d}+
N\nu^{-\alpha}\big(\dashint_{B^{+}_{\nu r}(z)}|D^{2}u|^{d}\,dx
\big)^{1/d},
\end{equation}
where $N$ depends only on $d$ and $\delta$.

\end{lemma}

The proof of this lemma coincides with that
of Lemma 2.5 of \cite{DKL} apart from the fact that
in place of Lemma 2.4 of \cite{Kr10} one should use
our Lemma \ref{lemma 10.5.2}. Also it is worth saying
that Lemma 2.3 of \cite{DKL} is contained in our Lemma
\ref{lemma 12.16.2}.

From now on we denote by $\alpha$ the smallest
of the $\alpha$'s in Lemmas \ref{lemma 10.5.1} and \ref{osca}.

The above arguments followed very closely the ones from
\cite{Kr10} and \cite{DKL}. At this point we 
will follow a different rout caused by the fact that
localization technique is not applicable
in our case as is pointed out in Section 
\ref{section 12.14.1}.

\begin{lemma}
                                       \label{lemma 10.9.1}
There exist constants
 $\rho_{0}\geq\rho_{1}>0$ depending only on $\Omega$
such that, for any $r>0$ and
 $\nu\geq64$ satisfying $\nu r\leq\rho_{1}$ and
 $z\in\Omega\setminus\Omega_{\rho_{0}}$ and
$u\in\WO^{2}_{d}(\Omega)$
we have
$$
\big(\dashint_{B_{r}(z)\cap\Omega}
\dashint_{B _{r}(z)\cap\Omega}
|D^{2}u(x)-D^{2}u(y)|^{\gamma}\,dxdy\big)^{1/\gamma}
$$
$$
\leq N\nu^{d/\gamma}\big(\dashint_{B_{ \nu r}(z)\cap\Omega}
(|F(D^{2}u)|^{d}+|Du |^{d})\,dx\big)^{1/d}
$$
\begin{equation}
                                             \label{10.9.1}
+N(\nu^{1+d/\gamma} r+\nu^{-\alpha})
\big(\dashint_{B_{ \nu r}(z)\cap\Omega}
|D^{2}u |^{d}\,dx\big)^{1/d},
\end{equation}
where  the constants
$N$ depend  only on $\Omega,d$, and $\delta$.

\end{lemma}

Proof. We take a $\rho_{0}>0$ for which at any point
$z_{0}\in\partial\Omega$ there is an orthonormal system
of coordinates  with the origin at $z_{0}$
 such that in the new coordinates $\tilde{x}
=(\tilde{x}_{1},\tilde{x}')$  there exists a function
$$
\psi \in C^{2}(\{\tilde{x}'\in\bR^{d-1}:
|\tilde{x}'|\leq4\rho_{0}\})
$$
with the $C^{2}$-norm controlled by a constant depending only
on $\Omega$ and such that
$$
\psi(0)=0,\quad \psi_{\tilde{x}_{i}}(0)=0,
\quad i=2,...,d,
$$
$$
 \{\tilde{x}:|\tilde{x}'|\leq4\rho_{0},
\psi(\tilde{x}')<\tilde{x}_{1}\leq\psi(\tilde{x}')+4\rho_{0}\}
\subset \Omega ,
$$
$$
\{\tilde{x}:|\tilde{x}'|\leq4\rho_{0},\tilde{x}_{1}
=\psi(\tilde{x}')\}\subset\partial\Omega.
$$
 
By decreasing $\rho_{0}$ if necessary we may assume that
for any 
$z$ with $\rho(z)\leq\rho_{0}$ there is a unique point
$z_{0}\in\partial\Omega$ such that $\rho(z)=|z-z_{0}|$.
Then in the system of coordinates associated with $z_{0}$
we have that $z=(|z-z_{0}|,0,...,0)$. 

We fix a $z$
with $\rho(z)\leq\rho_{0}$ and the above mentioned system
of new coordinates. Since $z$ is fixed and we are free
to use any orthonormal system of coordinates
we represent any point $x$ in $\bR^{d}$ as $x=(x_{1},x')$ 
and may assume that $z=(|z|,0,...0)$, $|z|\leq\rho_{0}$, and
 there exists a function
$$
\psi \in C^{2}(B_{4\rho_{0}}\cap\{x_{1}=0\})
$$
with the $C^{2}$-norm controlled by a constant depending only
on $\Omega$ and such that
$$
\psi(0)=0,\quad D_{i}\psi (0)=0,
\quad i=2,...,d,
$$
$$
\Gamma:=\{x:|x'|\leq4\rho_{0},
\psi(x')\leq x_{1}\leq\psi(x')+4\rho_{0}\}
\subset\bar{\Omega}
$$
$$
\{x:|x'|\leq4\rho_{0},x_{1}
=\psi(x')\}=\Gamma\cap\partial\Omega.
$$

Set 
$$
\hat{\Gamma}:=
\{y:|y'|<4\rho_{0},0\leq y_{1}\leq
4\rho_{0}\}
$$
and introduce a mapping $x\to y(x)$ of $\Gamma$ onto
$\hat{\Gamma}$ by
\begin{equation}
                                               \label{11.11.1}
x_{1}\to y_{1}=y_{1}(x)=x_{1}-\psi(x'),\quad
 x'\to y'=y'(x)=x'.
\end{equation}
This mapping has an inverse $y\to x(y)$. Since $D_{x'}\psi(0)=0$,
we can decrease $\rho_{0}$ if necessary, so that,
$$
|y(x^{1})-y(x^{2})|\leq 2|x^{1}-x^{2}|,\quad\forall x^{1},x^{2}
\in\Gamma,
$$
\begin{equation}
                                                         \label{12.4.1}
|x(y^{1})-x(y^{2})|\leq 2|y^{1}-y^{2}|,\quad\forall y^{1},y^{2}
\in\hat{\Gamma}.
\end{equation}

Next,  we take $\rho_{1}=\rho_{1}(\Omega)>0$
such that
$$
\rho_{1}\leq\rho_{0},\quad
B_{\rho_{1}}(z)\cap\Omega\subset \Gamma
$$
as long as $|z|\leq\rho_{0}$.
Observe that, owing to \eqref{12.4.1}, for $r\in(0,\rho_{1}]$ we have
\begin{equation}
                                             \label{11.11.2}
B_{ r/2}^{+}(|z|)\subset y\big(B_{r}(z)\cap\Omega
\big)\subset B_{2r}^{+}(|z|),
\end{equation}
where $B_{2r}^{+}(|z|)\subset\hat{\Gamma}$ since $2r\leq2\rho_{0}$
and $|z|+2r\leq4\rho_{0}$.
In terms of the inverse mapping $x=x(y)$ this is rewritten as
\begin{equation}
                                             \label{11.12.1}
x\big(B_{ r/2}^{+}(|z|)\big)\subset B_{r}(z)\cap\Omega
\subset x\big(B_{2r}^{+}(|z|)\big),\quad
r\leq\rho_{1}.
\end{equation}

Notice one more time that
$$
B_{r}^{+}(|z|)\subset \hat{\Gamma}
$$
for any $r\in(0,2\rho_{1}]$. We are going 
to use a few times the following consequence of 
\eqref{11.11.2} and \eqref{11.12.1}
and the fact that the Jacobian of the mapping $y=y(x)$
equals one:
\begin{equation}
                                               \label{11.12.2}
\dashint_{B^{+}_{ r/2}(|z|)}
|g(y)|\,dy\leq \dashint_{B_{r}(z)\cap\Omega}
|g(y(x))|\,dx
\leq  \dashint_{B^{+}_{2 r}(|z|)}
|g(y)|\,dy
\end{equation}
provided that $r\leq\rho_{1}$.

Finally, introduce a function
$$
\hat{u}(y)=u(x),\quad y=y(x), 
$$
which is well defined in the cylinder $\hat{\Gamma}$.

Now we apply Lemma \ref{osca} to 
 $\hat{u}(y)$. In order to avoid the confusion
while differentiating with respect to $y$ and $x$
we will supply the symbols of differentiation with
subscripts $y$ and $x$, respectively. Observe that
for any  $\nu\geq16$ and $r>0$
satisfying $ \nu r\leq 2\rho_{1}$
we have $B^{+}_{\nu r}(|z|))\subset\hat{\Gamma}$,
 which by Lemma \ref{osca} implies  that
$$
\big(\dashint_{B^{+}_{r}(|z|)}\dashint_{B^{+}_{r}(|z|)}
|D^{2}_{yy}\hat{u}(y^{1})-D^{2}_{yy}\hat{u}
(y^{2})|^{\gamma}\,dy^{1}dy^{2}\big)^{1/\gamma}
$$
\begin{equation}
                                              \label{10.10.2}
\leq N\nu^{d/\gamma}
\big(\dashint_{B^{+}_{\nu r}(|z|)}
|F(D^{2}_{yy}\hat{u})|^{d}\,dy\big)^{1/d}+
N\nu^{-\alpha}\big(\dashint_{B^{+}_{\nu r}(|z|)}|D^{2}_{yy}
\hat{u}|^{d}\,dy
\big)^{1/d}.
\end{equation}

Notice also  that for $y=y(x)$ and $x=x(y)$
$$
D_{y}\hat{u}(y)=(D_{x}u)(x)\frac{\partial x}{\partial y}(y),
$$
where $D_{y}$ and $D_{x}$ are row vectors and
$\partial x/\partial y$ is the matrix whose $ij$ entry is
$\partial x_{i}/\partial y_{j}$,
$$
D^{2}_{yy}\hat{u}(y)=[\frac{\partial x}{\partial y}(y)]^{*}
[D^{2}_{xx}u(x)]\frac{\partial x}{\partial y}(y)
+[D_{x_{k}}u(x)]D^{2}_{yy}x_{k}(y).
$$
Since $\partial x_{i}/\partial y_{j}(z)$ is the identity
matrix, for $y^{1},y^{2}\in B^{+}_{r}(|z|)$ we have
$$
|D^{2}_{yy}\hat{u}(y^{1})-D^{2}_{yy}\hat{u}
(y^{2})|\geq|D^{2}_{xx}u(x^{1})-D^{2}_{xx}u
(x^{2})|
$$
$$
-Nr (|D^{2}_{xx}u(x^{1})|+|D^{2}_{xx}u
(x^{2})|)-N(|Du(x^{1})|+|Du(x^{2})|),
$$
where $x^{i}=x(y^{i})$ and $N$ depends only on $\Omega$.
By using
\eqref{11.12.2} and observing that $r=2r/2$ and
$r/2\leq\rho_{1}$, we conclude that the left-hand side
of \eqref{10.10.2} is greater than or equal to
$$
\big(\dashint_{B_{ r/2}(z)\cap\Omega}
\dashint_{B_{ r/2}(z)\cap\Omega}|D^{2}_{xx}
u(x^{1})-D^{2}_{xx}u(x^{2})|^{\gamma}
\,dx^{1}dx^{2}\big)^{1/\gamma}
$$
$$
-N \dashint_{B_{2r}(z)\cap\Omega}
 (r|D^{2}_{xx}u|+|D_{x}u|)\,dx . 
$$
In what concerns the first term in the right-hand side of
\eqref{10.10.2} we have
$$
I_{r}(z):=\dashint_{B^{+}_{\nu r}(|z|)}
|F(D^{2}_{yy}\hat{u})|^{d}\,dy\leq N
\dashint_{B^{+}_{\nu r}(|z|)}
|F(D^{2}_{xx} u(x) )|^{d}\,dy
$$
$$
+N (\nu r)^{d}
\dashint_{B^{+}_{\nu r}(|z|)}|D^{2}_{xx}u(x)|^{d}\,dy
+N\dashint_{B^{+}_{\nu r}(|z|)}
|D _{x}u(x)|^{d}\,dy .
$$
By using \eqref{11.12.2} again  and the assuming that
$2\nu r\leq\rho_{1}$, we conclude that at point $z$
$$
I_{r}(z)
\leq N\dashint_{B_{2\nu r }(z)\cap\Omega}
(|F(D^{2}_{xx}u)|^{d}+(\nu r)^{d}|D^{2}_{xx}u |^{d}
+|D_{x}u |^{d})\,dx.
$$
One estimates the last term in \eqref{10.10.2} similarly
and concludes that for $\nu\geq16$ and $2\nu r\leq
 \rho_{1}$ it holds that
$$
\big(\dashint_{B_{ r/2}(z)\cap\Omega}
\dashint_{B_{ r/2}(z)\cap\Omega}|D^{2}_{xx}
u(x^{1})-D^{2}_{xx}u(x^{2})|^{\gamma}
\,dx^{1}dx^{2}\big)^{1/\gamma}
$$
$$
\leq N\nu^{d/\gamma}\big(\dashint_{B_{2\nu r }(z)\cap\Omega}
(|F(D_{xx}^{2}u)|^{d}+|D_{x}u |^{d})\,dx\big)^{1/d}
$$
$$
+N(\nu^{1+d/\gamma} r+\nu^{-\alpha})
\big(\dashint_{B_{2\nu r }(z)\cap\Omega}
|D^{2}_{xx}u |^{d}\,dx\big)^{1/d}.
$$
After that to obtain \eqref{10.9.1} it only remains to replace
$  r/2$ with $r$ and $4\nu$ with $\nu$. The lemma is proved.

\mysection{Global a priori estimates}

We take $ \rho_{1},\rho_{0}$ from Section \ref{section 10.7.2}
and do not assume that $F$ is independent of~$x$,
we only need it to satisfy  Assumption \ref{assump1}.
 
First we derive the following.

\begin{lemma}
                                         \label{lemma 10.6.4}
 For any $q\in[1,\infty)$ and $\mu>0$
there is a $\theta=\theta(d,\delta,K_{0},\mu,q)>0$ such that,
if Assumption \ref{assump1} is satisfied with this $\theta$,
then
for any $r\in(0,R_{0}]$ and $z\in\Omega$
$$
I(r,q,z):=\dashint_{\Omega\cap B_{r}(z)}\sup_{u''\in\cS,|u''|>t_{0}}
\frac{|F(u'',x)-\bar{F}(u'')|^{q}}{|u''|^{q}}\,dx\leq \mu^{q}.
$$
\end{lemma}

Proof. First observe that 
$|F(u'',x)|=|F(u'',x)-F(0,x)|\leq K_{0}|u''|$ and
$|\bar{F}(u'')|\leq \delta^{-1}|u''|$, so that 
$I(r,q,z)\leq N(\delta,K_{0},q)I(r,1,z)$ and we may assume that $q=1$.

Next, the functions $t^{-1}F(tu'',x)$
and $t^{-1}\bar{F}(t u'')$ are Lipschitz continuous
with respect to $u''$ with constants depending only
on $\delta$ and $K_{0}$.
Therefore, 
there exist points $u''(1),...,u''(n)$ with $n=n(\mu,d,\delta,K_{0})$,
 such that
$|u''(k)|=1$ and for any $u''\in\cS$ with $|u''|=1$ there exists  
a $k$ such that
$$
|t^{-1}F(tu'',x)-t^{-1}F(tu''(k),x)|\leq\mu/4,\quad
|t^{-1}\bar{F}(tu'' )-t^{-1}\bar{F}(tu''(k))|\leq\mu/4
$$
for all $t>0$.
Also note that  setting $t=|u''|$ we get
$$
\sup_{u''\in\cS,|u''|>t_{0}}
\frac{|F(u'',x)-\bar{F}(u'')|}{|u''|}=\sup_{u'':|u''|=1}
\sup_{t>t_{0}}
t^{-1}
|F(tu'' ,x)-\bar{F}(tu'')|
$$
$$
\leq\sum_{k=1}^{n}\sup_{t>t_{0}}
t^{-1}
|F(tu''(k) ,x)-\bar{F}(tu''(k))|+\mu/2.
$$
After that it is seen that our assertion is true with $q=1$ 
for
$\theta(d,\delta,K_{0},\mu,1)= \mu/(2n)$. The lemma is proved.

\begin{lemma}
                                       \label{lemma 11.27.5}
Let $r\in(0,\infty)$ and
 $\nu\geq 64$ satisfy  $\nu r\leq\rho_{1}\wedge
 R_{0} $.
Take   $\mu\in(0,\infty), \beta\in(1,\infty)$,
and suppose that Assumption \ref{assump1}
is satisfied with $\theta=\theta(d,\delta,K_{0},\mu,\beta d)$
(see Lemma \ref{lemma 10.6.4}).
Take a function $u\in\WO^{2}_{d}(\Omega)$  and denote
$$
I_{r}(z)=\big(\dashint_{B_{r}(z)\cap\Omega}
\dashint_{B _{r}(z)\cap\Omega}
|D^{2}u(x)-D^{2}u(y)|^{\gamma}\,dxdy\big)^{1/\gamma}.
$$

Then   for any $z\in\Omega$  
$$
I_{r} (z)\leq N\nu^{d/\gamma}\big(
\dashint_{B_{\nu r }(z)\cap\Omega}(|F[u]|^{d}
+|Du |^{d})\,dx\big)^{1/d}
$$
\begin{equation}
                                             \label{10.9.01}
+N(\mu\nu^{d/\gamma}+\nu^{1+d/\gamma} r+\nu^{-\alpha})
\big(
\dashint_{B_{\nu r }(z)\cap\Omega}
|D^{2}u |^{\beta' d}\,dx)^{1/(\beta'd)} +Nt_{0}\nu^{d/\gamma},
\end{equation}
where $\beta'=\beta/(\beta-1)$ and $N$ depends only on $\Omega,d,K_{0}$, 
and $\delta$.

\end{lemma}

Proof.   Take a $z\in\Omega$. If 
$z\in\Omega\setminus\Omega_{\nu r}$, then $z\in\Omega
\setminus\Omega_{\rho_{0}}$ since $\rho:=\nu r\leq \rho_{1}\leq
\rho_{0}$. Furthermore,   $\rho\leq R_{0}$. Therefore $\bar{F}=\bar{F}_{z,\rho}$ is well defined
 and  by Lemma \ref{lemma 10.9.1} we obtain
$$
I_{r}(z)
\leq N\nu^{d/\gamma}\big(\dashint_{B_{\rho}(z)\cap\Omega}
(|\bar{F}[ u]|^{d}+|Du |^{d})\,dx\big)^{1/d}
$$
\begin{equation}
                                                     \label{11.27.7}
+N(\nu^{1+d/\gamma} r+\nu^{-\alpha})
\big(\dashint_{B_{\rho}(z)\cap\Omega}
|D^{2}u |^{d}\,dx\big)^{1/d}.
\end{equation}
 Here 
$$
 \dashint_{B_{\rho}(z)\cap\Omega}
|\bar{F}[u]|^{d}\,dx \leq
N \dashint_{B_{\rho}(z)\cap\Omega}
|F[u]|^{d}\,dx +
N \dashint_{B_{\rho}(z)\cap\Omega}
|F[u]-\bar{F}[u]|^{d}\,dx ,
$$
where the last integral
is dominated by
$$
\dashint_{B_{\rho}(z) \cap\Omega}I_{|D^{2}u|>t_{0}}\frac{
|F[u]-\bar{F}[u]|^{d}}{|D^{2}u|^{d}}|D^{2}u|^{d}\,dx+Nt_{0}^{d},
$$
which in turn owing to Lemma \ref{lemma 10.6.4}
and H\"older's inequality is less than
$$
N\mu^{d}\big(\dashint_{B_{\rho}(z)\cap\Omega}
|D^{2}u|^{\beta'd}\,dx\big)^{1/\beta'}+Nt_{0}^{d}.
$$

It follows that
$$
\big(\dashint_{B_{\rho}(z)\cap\Omega}
|\bar{F}[u]|^{d}\,dx\big)^{1/d}\leq
N \big(\dashint_{B_{\rho}(z)\cap\Omega}
|F[u]|^{d}\,dx \big)^{1/d} 
$$
$$
+N\mu \big(\dashint_{B_{\rho}(z)\cap\Omega}
|D^{2}u|^{\beta'd}\,dx\big)^{1/(\beta'd)}+Nt_{0} .
$$
 
This and \eqref{11.27.7} yield \eqref{10.9.01} 
since 
$$
\big(\dashint_{B_{\rho}(z)\cap\Omega}
|D^{2}u |^{d}\,dx\big)^{1/d}\leq
\big(\dashint_{B_{\rho}(z)\cap\Omega}
|D^{2}u |^{\beta'd}\,dx\big)^{1/\beta'd}
$$
 by H\"older's inequality.

In case that $z\in\Omega_{\nu r}$ estimate 
\eqref{11.27.7} holds by Lemma \ref{lemma 10.5.2}
and  as above it leads to
\eqref{10.9.01}. The lemma is proved.

We now come to the main pointwise a priori estimate for 
nonlinear equations with VMO ``coefficients".
Introduce
\begin{equation}
                                                 \label{11.14.2}
h^{\#}_{\gamma }(z)=
\sup_{r>0}\big(\dashint_{B_{r}(z)\cap\Omega}
\dashint_{B_{r}(z)\cap\Omega}
|h(x)-h(y)|^{\gamma}\,dxdy\big)^{1/\gamma},
\end{equation}
$$
\bM h(z)=\sup_{r>0}\dashint_{B_{r}(z)\cap\Omega}
|h(x)|\,dx.
$$

\begin{theorem}
                                   \label{theorem 11.8.1}

Let $r\in(0,\infty)$ and
 $\nu\geq 64$ satisfy  $\nu r\leq\rho_{1}\wedge
  R_{0} $.
Take a $\mu\in(0,\infty), \beta\in(1,\infty)$,
and suppose that Assumption \ref{assump1}
is satisfied with $\theta=\theta(d,\delta,K_{0},\mu,\beta d)$.
Then for any function $u\in\WO^{2}_{d}(\Omega)$
we have in $\Omega$ that
$$
(D^{2}u)^{\#}_{\gamma } \leq N\nu^{d/\gamma}\bM^{1/d} (|F[u]|^{d}) 
+N\nu^{d/\gamma}\bM^{1/d} (|Du |^{d}) 
$$
$$
+N(\mu\nu^{d/\gamma}+\nu^{1+d/\gamma} r+\nu^{-\alpha})
\bM^{1/(\beta'd)} 
(|D^{2}u |^{\beta' d}) 
$$
\begin{equation}
                                             \label{11.27.8} 
+Nt_{0}\nu^{d/\gamma}+N r^{-d/\gamma} \|F[u]\|_{L_{d}(\Omega)},
\end{equation}
where $N$ depends only on $\Omega,d$, $K_{0}$, and $\delta$.

\end{theorem}

This theorem is an immediate consequence of
Lemma  \ref{lemma 11.27.5} and 
Corollary \ref{corollary 10.7.1}.
Indeed, the left-hand side of \eqref{11.27.8}  is the supremum
over $\rho>0$ of $I_{\rho}$. If $\rho\leq r$, $I_{\rho}$
is less than the right-hand side of \eqref{11.27.8}
by Lemma \ref{lemma 11.27.5}. However, if $\rho>r$,
then for $z\in\Omega$, obviously,
$$
I_{\rho}(z)\leq N\big(r^{-d}\int_{\Omega}
|D^{2}u|^{\gamma}\,dx\big)^{1/\gamma},
$$
which
is less than the right-hand side of \eqref{11.27.8}
by Corollary \ref{corollary 10.7.1} and the fact that
$F[u]=F[u]-F[0]=Lu$ for an $L\in\bL_{\delta,K_{0}}$.
 
Here is  the main a priori estimate.

\begin{theorem}
                                   \label{theorem 11.28.1}
Let $p\in(d,\infty)$. Then there exists
a constant $\theta>0$ depending only on $\Omega$, $p,d$, $K_{0}$,
 and $\delta$
such that if Assumption \ref{assump1}
is satisfied with this $\theta$,
then for any function
 $
u\in \WO^{2}_{p}(\Omega)
 $
we have
\begin{equation}
                                         \label{11.28.2}
\|u\|_{W^{2}_p(\Omega)}\le N
\| F[u]\|_{L_p(\Omega)} +Nt_{0},
\end{equation}
where $N$ depends only on  $\Omega,R_{0},d,p $, $K_{0}$, and $\delta$. 
 \end{theorem}

Proof. 
By Theorems \ref{theorem 11.14.1}, which is deferred to
Appendix, and Corollary
\ref{corollary 10.7.1}
 we have
$$
\|D^{2}u\|_{L_p(\Omega)}\leq N\|
(D^{2}u)^{\#}_{\gamma }\|_{L_{p}(\Omega)}+N\big(
\int_{\Omega}|F[u]|^{\gamma}\,dx\big)^{1/\gamma}
$$
$$
\leq N\|
(D^{2}u)^{\#}_{\gamma }\|_{L_{p}(\Omega)}+N\| F[u]\|_{L_p(\Omega)},
$$
where the last inequality follows from H\"older's inequality.

Then take $\beta\in(1,\infty)$ so that $\beta'd=(p+d)/2$
and take a $\mu>0$, which will be specified later,
and suppose that Assumption \ref{assump1}
is satisfied with $\theta=\theta(d,\delta,K_{0},\mu,\beta d)$
(see Lemma \ref{lemma 10.6.4}). Finally, take
$r\in(0,\infty)$ and
 $\nu\geq 64$ such that  $\nu r\leq\rho_{1}\wedge
 R_{0} $.

By Theorem   \ref{theorem 11.8.1}  and the Hardy-Littlewood
theorem
$$
\|
(D^{2}u)^{\#}_{\gamma }\|_{L_{p}(\Omega)}\leq N
(\nu^{d/\gamma}+r^{-d/\gamma})
 \|F[u]\|_{L_{p}(\Omega)}+N\nu^{d/\gamma}\|Du\|_{L_{p}(\Omega)} 
$$
$$
+N(\mu\nu^{d/\gamma}+\nu^{1+d/\gamma} r+\nu^{-\alpha})
\|D^{2}u\|_{L_{p}(\Omega)}+Nt_{0}\nu^{d/\gamma}.
$$
Hence,
$$
\|
D^{2}u\|_{L_{p}(\Omega)}\leq N
(\nu^{d/\gamma}+r^{-d/\gamma})
 \|F[u]\|_{L_{p}(\Omega)}+N\nu^{d/\gamma}\|Du\|_{L_{p}(\Omega)} 
$$
\begin{equation}
                                                    \label{11.28.5}
+N_{1}(\mu\nu^{d/\gamma}+\nu^{1+d/\gamma} r+\nu^{-\alpha})
\|D^{2}u\|_{L_{p}(\Omega)}+Nt_{0}\nu^{d/\gamma}.
\end{equation}
First we take and fix large $\nu\geq 64$ so that
$$
N_{1}\nu^{-\alpha}\leq1/4.
$$
Then we take and fix small $r>0$ so that $\nu r\leq\rho_{1}\wedge
  R_{0} $ and
$$
N_{1}\nu^{1+d/\gamma}r\leq1/4.
$$
Finally, we specify the value of $\mu>0$ we need so that
$$
N_{1}\mu\nu^{d/\gamma}\leq1/4.
$$
Then we conclude from \eqref{11.28.5} that
$$
\|
D^{2}u\|_{L_{p}(\Omega)}\leq N
 \|F[u]\|_{L_{p}(\Omega)}+N \|Du\|_{L_{p}(\Omega)} +Nt_{0}.
$$
After that to obtain \eqref{11.28.2} it only remains to use
interpolation inequalities and the Alexandrov
maximum principle, that says that $\max|u|\leq N\|F[u]\|
_{L_{p}(\Omega)}$. The theorem is proved.

\mysection{Proof of Theorem \protect\ref{theorem 10.5.1}}
                                               \label{section 11.29.1}

We only need to prove the  existence of solutions
and estimate \eqref{eq16.01}. We 
are going to use the following result of \cite{Kr13},
which is proved for any function $H$ satisfying
Assumption \ref{assumption 9.23.1} and such that
$$
\bar{H}:  =\sup_{x
\in\bR^{d}}|H(0,x)|<\infty. 
$$

Take a function $g\in C^{1,1}(\bar{\Omega})$.
\begin{theorem}
                                    \label{theorem 9.23.1}
 There are   constants $\hat{\delta}\in(0,\delta]$
and $\hat{K}_{0}\in[K_{0},\infty)$
depending only on $\delta$, $K_{0}$, and $d$ and there exists
a function $P(u) $ (independent of $x$), 
satisfying Assumption
\ref{assumption 9.23.1} with $\hat{\delta}$ and $\hat{K}_{0}$ in  place of 
$\delta$ and $K_{0}$
such that for any constant $K\geq0$  the equation
\begin{equation}
                                               \label{9.23.2}
\max(H[v],P[v]-K)=0
\end{equation}
in $\Omega$  (a.e.) with boundary condition $v=g$ on $\partial\Omega$
has a unique solution 
$v\in C^{0,1}(\bar{\Omega})\cap C^{1,1}_{loc}(\Omega)$.
In addition,  for all $i,j$, and $p\in(d,\infty)$,
\begin{equation}
                                                \label{1.13.1}
|v|,|D_{i}v|,\rho|D_{ij} v |\leq N(\bar{H}+K
+\|g\|_{C^{1,1}(\Omega)})\quad\text{in}
\quad \Omega \quad (a.e.), 
\end{equation}
\begin{equation}
                                                \label{1.13.2}
\|v\|_{W^{2}_{p}(\Omega)}\leq
N_{p}(\bar{H}+K+\|g\|_{W^{2}_{p}(\Omega)}),
\end{equation}
\begin{equation}
                                                \label{2.28.1} 
\|v\|_{C^{\alpha}(\Omega)}\leq N
(\|H[0]\|_{L_{d}(\Omega)}+\|g \|_{C^{\alpha}(\Omega)}),
\end{equation}
where 
$\alpha\in(0,1)$ is a constant depending
only on $d$ and $\delta$,
   $N$ is a constant depending only on $\Omega$ 
 and $\delta$, whereas $N_{p}$ only depends on the same objects and
$p$.

Finally, $P(u )$ is
constructed on the sole basis of $\delta$ 
and $d$, it is  positive homogeneous of degree one
and convex in $u$.

\end{theorem}

To derive Theorem \ref{theorem 10.5.1}  first assume that $g=0$,
 introduce a function of one variable by setting
$\xi_{K}(t)=0$ for $|t|\leq K$ and
$\xi_{K}(t)=t$ otherwise, set $H^{0}_{K}(x)=\xi_{K}(H(0,x))$, and define
$$
H_{K}(u,x)=\max(H(u,x)-H^{0}_{K}(x),P(u)-K),
$$

$$
F_{K}(u'',x)
=\max(F(u'',x),P(0,u'')-K),\quad
 G_{K}(u,x)=H_{K}(u,x)-F_{K}(u'',x).
$$
Notice that since $F(0,x)=0$, we have $|H^{0}_{K}|\leq \xi_{K}(\bar{G})
\leq\bar{G}$.
Also $|H_{K}(0,x)|\leq K$.

Below in this section by $N$ we denote
various constants which depend  only on $\Omega,R_{0},d,p $,
$K_{0}$, and $\delta$.
Observe that
$$
|G_{K}(u,x)|\leq |H(u,x)-H^{0}_{K}(x)-F(u'',x)|+|P(u)-P(0,u'')|
$$
$$
\leq K_{0}|u'|+2\bar{G}(x)
+\big|\sum_{k=0}^{d}u'_{k}\int_{0}^{1}P_{u'_{k}}(tu',u'')\,dt\big|
$$
\begin{equation}
                                                    \label{11.30.1}
\leq   N|u'|+2\bar{G}(x).
\end{equation}
Furthermore, $F_{K}$ obviously satisfies Assumption  
\ref{assump1} (iii) perhaps with  a constant $N$
(independent of $K$) in place of $K_{0}$.
 To check that the remaining conditions
in Assumption \ref{assump1} are satisfied take
$z\in\Omega$, $r\in(0,R_{0}]$, the function $\bar{F}=\bar{F}_{z,r}$
from Assumption \ref{assump1} and set
$$
\bar{F}_{K}(u'')=\max(\bar{F}(u''),P(0,u'')-K).
$$
Notice that
$$
t^{-1}
|F_{K}(tu'' ,x)-\bar{F}_{K}(tu'')| \leq
t^{-1}
|F(tu'' ,x)-\bar{F}(tu'')|,
$$
which implies that Assumption \ref{assump1} is satisfied
indeed with the same $\theta$.

By Theorem \ref{theorem 9.23.1} there is a unique
solution $v_{K}\in \WO^{2}_{p}(\Omega)$ of the equation
$$
H_{K}[v_{K}]=0.
$$
By Theorem \ref{theorem 11.28.1} 
\begin{equation}
                            \label{11.28.1}
\|v_{K}\|_{W^{2}_{p}(\Omega)}
\leq N \|F_{K}[v_{K}]\|_{L_{p}(\Omega)}+Nt_{0}
=N \|G_{K}[v_{K}]\|_{L_{p}(\Omega)}+Nt_{0}.
\end{equation}
It follows from \eqref{11.30.1} that
$$
\|G_{K}[v_{K}]\|_{L_{p}(\Omega)}\le 
N(\|Dv_{K}\|_{L_{p}(\Omega)}
+\|v_{K}\|_{L_{p}(\Omega)})+2\|\bar{G}\|_{L_{p}(\Omega)}.
$$
Therefore we obtain from \eqref{11.28.1}
that
\begin{equation}
                                                \label{11.28.3}
\|v_{K}\|_{W^{2}_{p}(\Omega)}
\leq N(\|Dv_{K}\|_{L_{p}(\Omega)}
+\|v_{K}\|_{L_{p}(\Omega)} +\|\bar{G}\|_{L_{p}(\Omega)}) 
+N t_{0} .
\end{equation}

Furthermore, $|H_{K}[0]|\leq |H[0]|\leq\bar{G}$ and there is
an operator $L\in\bL_{\hat\delta,\hat K_{0}}$ such that
$$
Lv_{K}=H_{K}[v_{K}]-H_{K}[0]=-H_{K}[0],
$$
which by the Alexandrov maximum principle implies that
$$
\|v_{K}\|_{L_{p}(\Omega)}\leq N\|\bar{G}\|_{L_{p}(\Omega)}.
$$
This and the interpolation inequalities allow us to
conclude from \eqref{11.28.3} that
\begin{equation}
                                                \label{11.28.4}
\|v_{K}\|_{W^{2}_{p}(\Omega)}
\leq N \|\bar{G}\|_{L_{p}(\Omega)} +N t_{0} .
\end{equation}
In this way we completed a crucial step consisting of
 obtaining a uniform
control of the $W^{2}_{p}(\Omega)$-norms of $v_{K}$.

We now let $K\to\infty$. As is well known, 
there is a sequence $K_{n}\to\infty$
as $n\to\infty$ and $v\in W^{2}_{p}(\Omega)$
  such that $v_{K}\to v$ weakly in $W^{2}_{p}(\Omega)$.
Of course, estimate \eqref{11.28.4} holds with $v$
in place of $v_{K}$, which yields \eqref{eq16.01}.

By the compactness of embedding of
$W^{2}_{p}(\Omega)$ into $C(\bar{\Omega})$ we have that
  $v_{K}\to v$ also uniformly.

Next, the operator $H[u]$ fits into the scheme of Section
5.6 of \cite{Kr85} since for any $u,v\in W^{2}_{p}(\Omega)$
there is an operator $L\in\bL_{\delta,K_{0}}$ such that
$H[u]-H[w]=L(u-w)$. Finally, by recalling that 
$|H^{0}_{K}|\leq \xi_{K}(\bar{G})$ we get
$$
|H[v_{K}]|=|\max(H[v_{K}]-H^{0}_{K},P[v_{K}]-K)-H[v_{K}]|
$$
$$
=|\max(0,P[v_{K}]-H[v_{K}]+H^{0}_{K}-K)-H^{0}_{K}|
$$
$$
\leq(P[v_{K}]-H[v_{K}]+H^{0}_{K}-K)_{+}+|H^{0}_{K}|
$$
$$
\leq(P[v_{K}]-H[v_{K}]+H^{0}_{K} -K)_{+}+ \xi_{K}(\bar{G})
$$
$$
\leq(N|D^{2}v_{K}|+
N|Dv_{K}|+N|v_{K}| +\bar{G}-K)_{+}+ \xi_{K}(\bar{G}),
$$
so that
$$
\| H[v_{K}] \|_{L_{d}(\Omega)}^{d}\leq NK^{d-p}
\int_{\Omega}
( |D^{2}v_{K}|+
 |Dv_{K}|+ |v_{K}|+\bar{G})^{p}\,dx\to0
$$
as $K\to\infty$. By combining all these facts and applying
  Theorems 3.5.15 and 3.5.6
of \cite{Kr85} we conclude that $H[v]=0$ and this finishes
the proof of   the theorem if $g\equiv0$. 

In the general case introduce
$$
\hat{H}(u,x)=H(u'_{0}+g(x),u'_{1}+D_{1}g(x),...,
u'_{d}+D_{d}g(x),u''_{ij}+D_{ij}g(x),x),
$$
$$
\hat{G}(u,x)=\hat{H}(u,x)-F(u'',x).
$$
Observe that
$$
|F[u+g]-F[u]|+|G[u+g]-G[u]|\leq N\bar{g},
$$
where
$$
\bar{g}=|D^{2}g|+|Dg|+|g|.
$$
It follows that
$$
|\hat{G}[u]|=|H[u+g]-F[u]|\leq|H[u+g]-F[u+g]|+
N\bar{g}
$$
$$
=|G[u+g]|+N\bar{g}\leq |G[u]|+N\bar{g}.
$$

We conclude that all our assumptions are satisfied for $\hat{H}$,
so that the equation $\hat{H}[u]=0$ has a unique solution $u\in
\WO^{2}_{p}(\Omega)$ and the corresponding estimate holds.
Then it only remains to set $v=u+g$. The theorem is proved.

 \mysection{Appendix}

In this section $\gamma$ is any number in $(0,1]$.

\begin{theorem}
                                       \label{theorem 11.14.1}
For any $p\in(1,\infty)$ 
 and $h\in L_{p}(\Omega)$
we have
\begin{equation}
                                                     \label{11.27.5}
\|h\|_{L_{p}(\Omega)}\leq N\|h^{\#}_{\gamma}\|_{L_{p}(\Omega)}
+N\big(\int_{\Omega}|h|^{\gamma}\,dx\big)^{1/\gamma},
\end{equation}
where $N$ depends only on $\gamma,d,p$, and $\Omega$.

\end{theorem}

Proof. It is convenient to supply the notation
$h^{\#}_{\gamma}$ with the subscript $\Omega$ reflecting the
fact that $h^{\#}_{\gamma}$ is defined by \eqref{11.14.2}
for each particular $\Omega$. Therefore, in this proof
we denote the right-hand side of \eqref{11.14.2} by
$h^{\#}_{\Omega,\gamma}$. The notation $\bM_{\bR^{d}}$
has a similar meaning.

Next, it is well known that on account of $\Omega$ being a
bounded domain of class $C^{2}$ there is a $\rho_{0}>0$
depending only on $\Omega$ such that in
$\Omega^{2\rho_{0}}\setminus\Omega$, where
$$
\Omega^{2\rho_{0}}=\{x :\text{dist}\,(x, 
\Omega)<2\rho_{0}\},
$$
there is a $C^{2}$ mapping,
 which maps $\Omega^{2\rho_{0}}\setminus\Omega$
onto $\bar\Omega\setminus\bar\Omega_{2\rho_{0}}$ in a one-to-one
way with a $C^{2}$ inverse and preserves
$\partial\Omega$. 
We continue this mapping inside $\Omega$ as the identity
mapping and call   $\psi(x)$ the such obtained mapping of 
$\Omega^{2\rho_{0}} $ into itself. Of course,
$\psi$ will not be of class $C^{2}$, but yet its Lipschitz
constant is finite and depends only on $\Omega$.
Then take a $\zeta\in
C^{\infty}_{0}(\Omega^{\rho_{0} } )$ such that
$\zeta=1$ on $\Omega$, $0\leq\zeta\leq1$, and define
$$
\hat{h}(x)=h(\psi(x))\zeta^{2}(x),\quad x\in\Omega^{2\rho_{0}}
$$ 
(and continue $\hat{h}$ as zero outside $
\Omega^{2\rho_{0}} $ where formally $\psi(x)$
is not defined).

Now take a $z\in\bR^{d}$ and assume that
\begin{equation}
                                             \label{11.15.2}
d(z):=\dist(z,\Omega)\geq2\rho_{0}.
\end{equation}
If
$r\in(0, d(z) -\rho_{0}]$, then
$$
I_{r}(z):=\big(\dashint_{B_{r}(z)}\dashint_{B_{r}(z)}|
\hat{h}(x)-\hat{h} (y)|^{\gamma}\,dxdy\big)^{1/\gamma}=0.
$$
In case \eqref{11.15.2} holds and
$r>d(z) -\rho_{0}$ we have $r^{-d}\leq  
(d(z) -\rho_{0})^{-d}\leq N(1+|z|)^{-d}$ and
$$
I_{r}(z)\leq Nr^{-d/\gamma} 
\big(\int_{\bR^{d}}|\hat{h}|^{\gamma}\,dx\big)^{1/\gamma}
\leq N(1+|z|)^{-d/\gamma}
\big(\int_{\Omega}|h|^{\gamma}\,dx\big)^{1/\gamma}.
$$
Generally, if $r\geq \rho>0$, then
\begin{equation}
                                              \label{11.27.4}
I_{r}(z)\leq N\rho^{-d/\gamma}
\big(\int_{\Omega}|h|^{\gamma}\,dx\big)^{1/\gamma}.
\end{equation}

Next assume that
\begin{equation}
                                             \label{11.27.1}
\rho(z)<2\rho_{0}.
\end{equation}
In that case observe that
$$
|h(\psi(x))\zeta^{2}(x)-h(\psi(y))\zeta^{2} (y)|\leq
|h(\psi(x))I_{\Omega^{2\rho_{0}}}(x)|\,|\zeta^{2}(x)-\zeta^{2} (y)|
$$
$$
+|h(\psi(x))I_{\Omega^{2\rho_{0}}}(x)- 
h(\psi(y))I_{\Omega^{2\rho_{0}}}(y)|\zeta^{2} (y),
$$
where the first term is estimated by 
$$
N|x-y|\,|h(\psi(x))I_{\Omega^{2\rho_{0}}}(x)|
$$
and the second one equals
$$
|h(\psi(x))I_{\Omega^{2\rho_{0}}}(x)- 
h(\psi(y))I_{\Omega^{2\rho_{0}}}(y)|\zeta  (y)[\zeta(y)
-\zeta(x)]
$$
$$
+|h(\psi(x))I_{\Omega^{2\rho_{0}}}(x)- 
h(\psi(y))I_{\Omega^{2\rho_{0}}}(y)|\zeta  (y)\zeta  (x).
$$
It follows that
$$
I^{\gamma}_{r}(z)\leq
Nr^{\gamma}
\bM_{\bR^{d}}(h^{\gamma}(\psi)I_{\Omega^{2\rho_{0}}})(z)+
 Nr^{-2d}J_{r}(z),
$$
where
$$
J_{r}(z)=\int_{B_{r}(z)\cap\Omega^{2\rho_{0}}}
\int_{B_{r}(z)\cap\Omega^{2\rho_{0}}}
|h(\psi(x))-h(\psi(y))|^{\gamma}\,dxdy .
$$
We fix $z$ and $r$ and represent $J_{r}(z)$ as 
$$
J_{r}(z)=\sum_{i,j=1}^{2}J^{ij} 
$$
 according to integrating
with respect to $x\in \Gamma_{i}$ and $y\in\Gamma_{j}$,
where
$$
\Gamma_{1}=B_{r}(z)\cap(\Omega^{2\rho_{0}}\setminus\Omega),
\quad \Gamma_{2}=B_{r}(z)\cap \Omega .
$$
Notice that if $J_{r}^{22}\ne0$, then
 $\Gamma_{2}\ne\emptyset$,  $|z-\psi(z)|
\leq Nr$, 
\begin{equation}
                                                    \label{12.15.1}
\Gamma_{2}\subset B_{Nr}(\psi(z))\cap \Omega ,
\end{equation}  
and since
also
$$ 
|B_{Nr}(\psi(z))\cap \Omega|\leq Nr^{d},
$$  we have
\begin{equation}
                                              \label{11.27.3}
(r^{-2d}J_{r}^{ij})^{1/\gamma}\leq Nh^{\#}_{\Omega,\gamma}
(\psi(z))
\end{equation}
for $i=j=2$.
If $J_{r}^{12}\ne0$, then we have \eqref{12.15.1} and also
 $$
\psi(\Gamma_{1})\in B_{Nr}(\psi(z))\cap \Omega .
$$
By changing variables of integration we see that
\eqref{11.27.3} holds for $i\ne j$ as well.
Finally, regardless $J_{r}^{11}=0$ or not, changing 
variables shows that \eqref{11.27.3} holds in the
remaining case when $i=j=1$.

Hence, for $z$ satisfying \eqref{11.27.1} we have
$$
I_{r}(z)\leq N h^{\#}_{\Omega,\gamma}
(\psi(z))+Nr
\bM^{1/\gamma}_{\bR^{d}}(h^{\gamma}(\psi)I_{\Omega^{2\rho_{0}}})(z),
$$
which along with \eqref{11.27.4} shows that for any $\rho>0$
in both cases: $r\leq\rho$ and $r>\rho$, we have
$$
I_{r}(z)\leq N h^{\#}_{\Omega,\gamma}
(\psi(z))+N\rho
\bM^{1/\gamma}_{\bR^{d}}(h^{\gamma}(\psi)I_{\Omega^{2\rho_{0}}})(z)
+N\rho^{-d/\gamma}\big(\int_{\Omega}|h|^{\gamma}\,dx\big)^{1/\gamma},
$$
provided that  $z$ satisfies \eqref{11.27.1}.

This and \eqref{11.27.4} imply that in any case for any $\rho>0$
$$
\hat{h}^{\#}_{\bR^{d},\gamma}(z)\leq N h^{\#}_{\Omega,\gamma}
(\psi(z))+N\rho
\bM^{1/\gamma}_{\bR^{d}}(h^{\gamma}(\psi)I_{\Omega^{2\rho_{0}}})(z)
$$
\begin{equation}
                                             \label{11.15.1}
+N(1+\rho^{-d/\gamma})(1+|z|)^{-d/\gamma}
\big(\int_{\Omega}|h|^{\gamma}\,dx\big)^{1/\gamma},
\end{equation}
where $N$ depends only on $\Omega$.

Due to Theorem 5.3 of \cite{Kr10}
$$
\|\hat{h}\|_{L_{p}(\bR^{d})}\leq N\|\hat{h}^{\#}_{\bR^{d},\gamma}
\|_{L_{p}(\bR^{d})},
$$
where $N$ depends only on $p,\gamma,d$, which combined with the
  Hardy-Littlewood
theorem and \eqref{11.15.1} implies that
$$
\|h\|_{L_{p}(\Omega)}\leq N\|h^{\#}_{\Omega,\gamma}\|_{L_{p}(\Omega)}
+N(1+\rho^{-d/\gamma})
\big(\int_{\Omega}|h|^{\gamma}\,dx\big)^{1/\gamma}
+N_{1}\rho\|h\|_{L_{p}(\Omega)},
$$
where the constants depend only on $\Omega,\gamma$, and $p$.
Choosing $\rho$ so that $N_{1}\rho=1/2$ yields
\eqref{11.27.5} and the theorem is proved.

\end{document}